\documentclass[aip,graphicx]{revtex4-1}

\draft 

\usepackage{graphicx}
\usepackage{dcolumn}
\usepackage{bm}

\usepackage[utf8]{inputenc}
\usepackage[T1]{fontenc}
\usepackage{mathptmx}

\usepackage{amssymb}

\newcommand{\tr}{tr}
\newcommand{\id}{id}
\newcommand{\arctanh}{$arctanh$ }

\begin{document}

\newtheorem{lm}{Lemma}
\newtheorem{theorem}{Theorem}
\newtheorem{df}{Definition}
\newtheorem{prop}{Proposition}
\newtheorem{remark}{Remark}


\title{Global and local bifurcations, three-dimensional Henon maps and discrete Lorenz attractors}

\author{I. Ovsyannikov}
\email[]{iovsyann@uni-bremen.de}
\email[]{ivan.i.ovsyannikov@gmail.com}
 \affiliation{University of Bremen, MARUM, Department of Mathematics, Germany.}
 
 \affiliation{Lobachevsky State University of Nizhny Novgorod, ITMM, Russia}



\date{\today}

\begin{abstract}
Lorenz attractors play an important role in the modern theory of dynamical systems. The reason is that they are robust, i.e.  preserve their chaotic properties under various kinds of perturbations. This means that such attractors can exist in applied models and be observed in experiments. It is known that discrete Lorenz attractors can appear in local and global bifurcations of multidimensional diffeomorphisms. However, to date only partial cases were investigated.
In this paper bifurcations of homoclinic and heteroclinic cycles with quadratic tangencies of invariant manifolds are studied. A full list of such bifurcations, leading to the appearance of discrete Lorenz attractors is provided. In addition, with help of numerical techniques, it was proved that if one reverses time in the diffeomorphisms described above, the resulting systems also have such attractors. This result is an important step in the systematic studies of chaos and hyperchaos.
\end{abstract}

\maketitle

\begin{quotation}
Dynamical chaos is a topic that attracts a high interest of researchers as
dynamical models with complex behaviour can be widely met in applications. This relates, for example, to climate models \cite{L86, P98}, 
fluid mechanics \cite{PF80, P11, EM18}, nonholonomic rigid body mechanics \cite{GG15}, laser dynamics \cite{VV93} and many others.
Usually, the presence of chaos
in a dynamical system is connected to the existence of strange attractors. According to Afraimovich and Shilnikov \cite{AS83}, they
can be divided into two main types -- {\em genuine strange attractors} and {\em quasiattractors}. Many of well-known
types of chaotic attractors, such as H\'enon-like attractors, spiral attractors, R\"ossler attractors, attractors in the
Chua circuits etc., are quasiattractors in the sense that arbitrary small perturbations of them lead to appearance of stable periodic orbits. 
This is not the case for genuine strange attractors,
they exist in open domains in the space of dynamical systems and even in the case when they are structurally unstable (e.g. the Lorenz attractor),
stable periodic orbits are not born in bifurcations.
\end{quotation}

\section{\label{sec:intro}Introduction}

Among the known genuine strange attractors there are hyperbolic attractors, Lorenz attractors and wild hyperbolic attractors. The latter are remarkable
by the fact that they, unlike the previous ones, allow homoclinic tangencies and therefore they belong to
Newhouse domains \cite{NPT, GST93b}.
In addition, such attractors are stable, closed and chain-transitive invariant sets. Chain-transitivity means that any point of attractor $\Lambda$ 
is admissible by $\varepsilon$-orbits from any other point of $\Lambda$; stability means the existence of an open absorbing domain containing the attractor such that any orbit 
	entering the domain tends to $\Lambda$ exponentially fast.
All above constitutes the definition of an attractor by Ruelle and Conley \cite{R81, C78}.

An example of a wild spiral (hyperbolic) attractor was presented first by Turaev and Shilnikov \cite{TS98}.
Another important example of a wild hyperbolic attractor is the discrete
Lorenz attractor which appears, in particular, in the Poincare maps for periodically perturbed flows with Lorenz attractors \cite{TS08}.
They were also observed in applications, such as the nonholonomic rattleback model \cite{GG15} and two-component convection \cite{EM18}.
It is well-known that the classical Lorenz attractor does not allow homoclinic tangencies \cite{ABS77, ABS83}, however the latter can
appear under small non-autonomous periodic perturbations, this is a method to make a wild hyperbolic attractor out of the Lorenz attractor. 
The reason why stable periodic orbits do not arise from  bifurcations of these tangencies, is that Lorenz-like attractors possess a 
pseudo-hyperbolic structure, and this property is preserved under small perturbations. Recall that the pseudo-hyperbolic structure exists in map $f$ if its differential $Df$  in the restriction onto the absorbing domain ${\cal D}$ of attractor $\Lambda$ admits  for any point $x\in {\cal D}$ an invariant splitting of form 
$E^{ss}_x\oplus E^{uc}_x$, such that $Df$ is strongly contracting along directions $E^{ss}$ and expands volumes in transversal to $E^{ss}$ sections $E^{uc}$ (see Refs. \onlinecite{TS98, TS08} for details), and this splitting continuously depend on the point $x$.

One of the peculiarities of discrete Lorenz attractors is that such attractors can be born at local bifurcations
of periodic orbits having three or more
multipliers lying on the unit circle. Thus the corresponding attractors can be found in particular
models which have a sufficient number of parameters to provide the mentioned degeneracy. The following 3D H\'enon map
\begin{equation}
	\bar x = y, \;\; \bar y = z, \;\; \bar z = M_1 + B x + M_2 y - z^2
	\label{H3D}
\end{equation}
which depends on three  parameters $M_1$, $M_2$ and $B$ and has constant Jacobian $B$, is an example of such a model. 
In Refs. \onlinecite{GOST05, GMO06, GGOT13}
it was shown that  map (\ref{H3D}) possesses a discrete Lorenz attractor in some open parameter domain near 
point $(M_1, B, M_2) = (1/4,1,1)$, where the map has a fixed point with the triplet $(-1,-1,+1)$ of multipliers. Recently in Ref.~\onlinecite{GKKST}, it was proved that map~(\ref{H3D}) also has discrete Lorenz attractors in the orientation reversing case $B < 0$, near the codimension-three bifurcation $(M_1, M_2, B) = (7/4, -1, -1)$, when the map has a fixed point with eigenvalues $(i, -i, -1)$. The properties of these attractors were studied in Ref.~\onlinecite{GGKS21}. Note that the second iterate of the map in this case has a fixed point with multipliers $(-1,-1,+1)$, but, unlike the orientable case, they form three Jordan blocks instead of two, so the bifurcations that happen for $B > 0$ and $B < 0$ are principally different.

These results immediately imply the birth of discrete Lorenz attractors in systems where map (\ref{H3D}) appears, e.g. as a Poincare map. This, in particular, happens in global (homoclinic and heteroclinic) bifurcations. 
The first such example was considered in Ref.~\onlinecite{GMO06} for a homoclinic tangency to a saddle-focus. 
Later analogous results 
were obtained for heteroclinic cycles containing saddle-foci \cite{GST09, GO10, GO13} and homoclinic and heteroclinic cycles consisting of saddles and having additional degeneracies, such that non-simple homoclinic (heteroclinic) orbits~\cite{GOT14, IO21}, or a resonance condition on the eigenvalues at the fixed point~\cite{GO17}. Note that the presence of 
saddle-foci or degeneracies of certain kinds in these cases is a very important condition for the birth of Lorenz-like attractors as it prevents 
from the existence of lower-dimensional center manifolds and makes the dynamics to be effectively three-dimensional (see Refs. \onlinecite{T96, GST93c, Tat01}).
Another important condition was imposed on the values of the Jacobian of the map evaluated in the fixed points. It is based on the fact that
the orbits under consideration may spend unboundedly large number of iterations in the neighbourhoods of saddle fixed points. In the homoclinic case
this means that if the Jacobian is separated from unity, the phase volumes near such orbits will be either unboundedly expanded or unboundedly contracted,
and the dynamics will have effective dimension less than three. In the same way, for the heteroclinic cases it is necessary to demand that all
the Jacobians are not simultaneously contracting ($< 1$) or simultaneously expanding ($> 1$). Thus, in order to get Lorenz attractors in bifurcations of heteroclinic
cycles, one needs to consider ``contracting-expanding'' or ``mixed'' cases.

The results mentioned above concern homoclinic and heteroclinic $(2,1)$--cycles, i.e. those containing fixed points with a two-dimensional stable manifold and one-dimensional unstable manifold. In a three-dimensional space two other kinds of cycles are possible:  $(1,2)$--cycles, consisting of fixed points of type $(1,2)$, and heterodimensional cycles, having saddles of both stability types. Heterodimensional cycles are out of scope of this paper. Every system with a $(1,2)$--cycle can be regarded as an inverse to a system with a $(2,1)$--cycle, in which the first return map is the 3D Henon map~(\ref{H3D}). This immediately means that the first return map in $(1,2)$--cycles will be close up to small terms to the inverse map of~(\ref{H3D}), this is the following three-dimensional Henon map:
\begin{equation}\label{Henon2}
	\bar x = y, \;\;
	\bar y = z, \;\;
	\bar z = \hat M_1 + \hat B x + \hat M_2 z - y^2.
\end{equation}
This map was  obtained first in Ref.~\onlinecite{GST93c} as a first return map along a homoclinic orbit to a four-dimensional saddle-focus of $(2,2)$ type. As map (\ref{Henon2}) is an inverse to (\ref{H3D}), there exist domains in the space of parameters $(\hat M_1, \hat M_2, \hat B)$, in which it possesses a discrete Lorenz repeller \cite{GOST05, GKKST, GGKS21}. However, Lorenz attractors have not been found in this map before.

In the present paper, bifurcations, leading to the birth of discrete Lorenz attractors in homoclinic and heteroclinic cycles of type $(2,1)$ and $(1,2)$ with quadratic tangencies of invariant manifolds, are studied.
The full list of such bifurcations is presented, namely homoclinic (or heteroclinic, consisting of two fixed points) cycles, having
\begin{itemize}
	\item saddle-focus fixed points;
	\item non-simple homoclinic (heteroclinic) orbits;
	\item resonant fixed points.
\end{itemize} 
This list includes 
known results for orientable maps~\cite{GMO06, GST09, GO10, GO13, GOT14, IO21, GO17}, which were extended here to non-orientable maps. The latter is possible due to the recent result \cite{GKKST} that map (\ref{H3D}) possesses discrete Lorenz attractors also in the non-orientable case $B < 0$.
 
In addition, new cases are considered, when the cycle contains resonant fixed points that undergo the Belyakov transition from saddle to saddle-focus. This happens when at the bifurcation moment the stable multipliers have multiplicity two, and under small perturbations such a pair splits either in two different real eigenvalues or in a complex-conjugate pair.

The results of the paper are the following. First, for $(2,1)$--cycles it is shown that in the space of dynamical systems the original system is a limit of a sequence of open subsets containing systems with discrete Lorenz attractors, Theorem~\ref{thmmain}. The proof is based on the fact (Lemma~\ref{resc_lemma}) that the first return map near a homoclinic or heteroclinic cycle can be represented in the form of three-dimensional Henon map (\ref{H3D}), which has the discrete Lorenz attractor, see Refs.~\onlinecite{GOST05, GMO06, GGOT13, GKKST}. Then, it is proved numerically (Lemma~\ref{Lemma_H2}), that the inverse to (\ref{H3D}) map~(\ref{Henon2}) possesses the discrete Lorenz attractor near certain period-$6$ points. This result implies that systems with $(1, 2)$--cycles are also limits of sequences of open subsets in which systems have discrete Lorenz attractors

The paper is organised as follows. Section~\ref{sec:def} contains the statement of the problem, main definitions, and also the main results, Theorems~\ref{thmmain} and \ref{conj_inv} are formulated there. In Section~\ref{sec:Henon} the birth of discrete Lorenz attractors in the inverse 3D Henon map (\ref{Henon2}) is studied. In Section~\ref{sec:fret} the first return map is constructed. For all cases under consideration, local and global maps are written. At the end of the Section, there is the rescaling Lemma~\ref{resc_lemma}, stating that the first return map for all $(2,1)$--cycles coinsides with the 3D Henon map (\ref{H3D}).
\begin{figure}
\includegraphics[width=7cm]{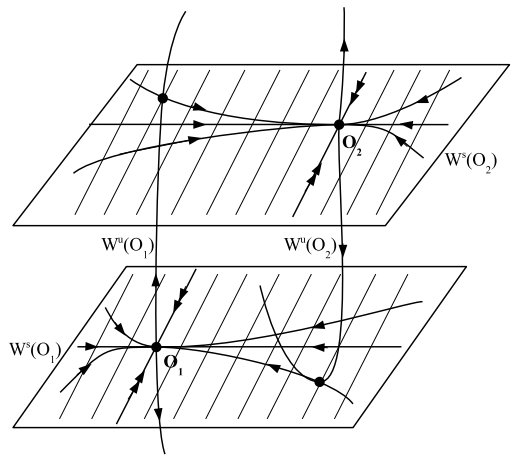}
\caption{\label{saddles} A heteroclinic cycle consisting of two saddles, with a quadratic tangency of manifolds.} 
\end{figure}



\section{Statement of the problem and main results}\label{sec:def}

In this section two classes of problems are set up: 
the study of bifurcations in homoclinic and heteroclinic cycles. These two classes are defined by corresponding conditions {\bf A}--{\bf C}, and fall into subcases, defined by an additional condition {\bf D}.

Let a $C^r$-smooth three-dimensional diffeomorphism $f_0$, $r \geq 4$, satisfy the following conditions:\\

{\bf I. The homoclinic case.}

\begin{itemize}
\item[\bf A.] $f_0$ has a fixed point $O$ with multipliers 
$(\lambda_1, \lambda_2, \gamma)$, where $|\lambda_{1,2}| < 1 < |\gamma|$;

\item[\bf B.] $|J(O)| \equiv |\lambda_1 \lambda_2 \gamma| = 1$;

\item[\bf C.] Invariant manifolds $W^u(O)$ and $W^s(O)$ have a quadratic tangency at the points of a homoclinic orbit $\Gamma_0$.
\end{itemize}

{\bf II. The heteroclinic case.}

\begin{itemize}
\item[\bf A.] $f_0$ has two fixed points
$O_1$ and $O_2$ of type $(2,1)$, i.e. each $O_j$ has multipliers $(\lambda_{(j)1}, \lambda_{(j)2}, \gamma_{(j)})$ with
$|\lambda_{(j)1, (j)2}| < 1 < |\gamma_{(j)}|$ for $j = 1,2$.

\item[\bf B.] The absolute value of the Jacobian of $f_0$ is less than one in one fixed point and greater than one in another one. Without loss of generality, $|J(O_1)| \equiv |\lambda_{(1)1} \lambda_{(1)2} \gamma_{(1)}| < 1$ and  $|J(O_2)| \equiv |\lambda_{(2)1} \lambda_{(2)2} \gamma_{(2)}| > 1$.

\item[\bf C.] There exists a heteroclinic cycle such that one-dimensional stable manifolds $W^u(O_1)$
intersect transversely two-dimensional stable manifolds $W^s(O_2)$
in the points of hereroclinic orbits $\Gamma_{12}$ and 
unstable manifold $W^u(O_2)$ has a quadratic tangency with stable manifold $W^s(O_1)$ 
at the points of a non-transversal heteroclinic orbit $\Gamma_{21}$.
\end{itemize}


The goal of this paper is the study of bifurcations of single-round periodic orbits, lying in a some small neighborhood of the
homoclinic or the heteroclinic cycle defined by conditions 
{\bf A}--{\bf C}. The main attention is paid to codimension-three bifurcations leading to the appearance of discrete Lorenz attractors. For this purpose, the first return map $T$ along the cycle is constructed, such that  single-round periodic orbits become  fixed points of map $T$. This map is three-dimensional, however, in some circumstances it can have lower-dimensional invariant submanifolds, that prevent from the existence of Lorenz-like attractors. To avoid this, 
the effective dimension of the problem~\cite{T96} should be kept equal to three. This is achieved by imposing an additional
condition {\bf D}, which gives the following subcases.  
\\

{\bf D.I. The homoclinic case}
\begin{itemize}
\item[\bf 1.] Point $O$ is a saddle-focus, i.e. $\lambda_{1,2} = \lambda e^{\pm i \varphi}$, $0 < \varphi < \pi$, see Ref.~\onlinecite{GMO06};
\item[\bf 2.] The quadratic tangency at $\Gamma_0$ is non-simple (the definitions are given below), see Ref.~\onlinecite{GOT14}:
  \begin{itemize}
  \item[\bf 2.a.] the homoclinic orbit undergoes an inclination flip;
  \item[\bf 2.b.] the homoclinic orbit undergoes an orbit flip;
  \end{itemize}
\item[\bf 3.] The multipliers at $O$ satisfy the following conditions:
  \begin{itemize}
  \item[\bf 3.a.] an alternating resonance $\lambda_1 = -\lambda_2 = \lambda$, see Ref.~\onlinecite{GO17};
  \item[\bf 3.b.] a Belyakov transition from saddle to saddle-focus $\lambda_1 = \lambda_2 = \lambda$.
  \end{itemize}
\end{itemize}

{\bf D.II. The heteroclinic case}
\begin{itemize}
\item[\bf 1.] One of the points $O_1$ and $O_2$, or they both, are saddle-foci, see Refs.~\onlinecite{GST09, GO10, GO13};
\item[\bf 2.] One of the heteroclinic orbits is non-simple, see Ref.~\onlinecite{IO21}:
  \begin{itemize}
  \item[\bf 2.a.] the non-transversal orbit $\Gamma_{21}$ undergoes an inclination flip;
  \item[\bf 2.b.] the non-transversal orbit $\Gamma_{21}$ undergoes an orbit flip;
  \item[\bf 2.c.] the transversal orbit $\Gamma_{12}$ undergoes an orbit flip;
  \end{itemize}
\item[\bf 3.] The stable multipliers at one fixed point $O_j$, $j = 1$ or $2$, satisfy the following condition
  \begin{itemize}
  \item[\bf 3.a.] an alternating resonance $\lambda_{(j)1} = -\lambda_{(j)2} = \lambda$;
  \item[\bf 3.b.] a Belyakov transition from saddle to saddle-focus $\lambda_{(j)1} = \lambda_{(j)2} = \lambda$.
  \end{itemize}
\end{itemize}

The notation of the cases is the following: the Roman numbers {\bf I} and {\bf II} denote the homoclinic and heteroclinic cases respectively, and arabic numbers with letters denote the subcases given by condition {\bf D}. For example, case {\bf I.2.a} is a homoclinic tangency with an inclination flip to a saddle fixed point, and case {\bf II.1} is
a heteroclinic cycle consisting of a saddle and a saddle-focus or two saddle-foci. All the required definitions are given later in this section.

The following lemma states that fulfilment of at least one subcase of {\bf D} is a necessary condition for discrete Lorenz attractors to appear in bifurcations. 

\begin{lm}\label{lem_notD}
Assume that $f_0$ satisfies conditions {\bf A}--{\bf C}, but condition {\bf D} is not fulfilled. Then in small generic unfoldings $f_\mu$ of $f_0$,
the first return map along the homoclinic or heteroclinic cycle possesses a global invariant manifold of dimension lower than three.
\end{lm}

Diffeomorphisms close to $f_0$ and satisfying conditions {\bf A}--{\bf D} compose locally connected bifurcation surfaces in the space of $C^r$-diffeomorphisms. They have codimension one
in case {\bf II.1} codimension two in cases {\bf I.1}, {\bf II.2--3} and codimension three in cases {\bf I.2--3}. 

Generic bifurcations of systems close to $f_0$ are studied in three-parametric
unfoldings, that are, the families of diffeomorphisms 
$f_\mu$, $\mu = (\mu_1, \mu_2, \mu_3)$ such that $\left. f_\mu \right|_{\mu = 0} = f_0$.
The first parameter $\mu_1$ is selected as the splitting distance of the quadratic tangency (homoclinic or heteroclinic), defined by condition {\bf C}. The second parameter $\mu_2$ controls condition {\bf D} such that in the
cases {\bf I--II.1}, when one or more saddle-foci are present, $\mu_2$ is monotonically related to the complex argument $\varphi_j$ of stable multipliers of one saddle-focus $O_j$: 
$\mu_2 = \varphi_j - \varphi_0$, and in all other cases $\mu_2$ unfolds the corresponding degeneracy given by condition {\bf D}.

The third parameter $\mu_3$ controls the Jacobians at the saddle points, such that in the homoclinic case {\bf I}
it is responsible to a deviation of the Jacobian from unity: 
\begin{equation} \label{mu3_hom}
\mu_3 = |J(O)| - 1,
\end{equation}%
and in the heteroclinic case~{\bf II} it is a deviation of a functional $S$:
\begin{equation} \label{mu3_het}
	\mu_3 = S(f_\mu) - S(f_0),
\end{equation}
where $$S = -\frac {\ln |J(O_{1})|}{\ln |J(O_{2})|}.$$ 

Note that when the system has a homoclinic or heteroclinic cycle containing saddle-foci (cases {\bf I--II.1}), 
there exist continuous invariants of topological conjugacy on the set of non-wandering orbits  
($\Omega$-moduli), and the most important moduli are the angular arguments \cite{G00a} of complex eigenvalues. 
This is the reason why at least one of them should be chosen as a control parameter $\mu_2$. 
The condition {\bf B} on Jacobians are essential for all cases, as the orbits
lying in a small neighbourhood of the homoclinic or heteroclinic cycle, 
may spend arbitrary large number of iterations near saddle points, and if the absolute value of the 
Jacobian in the homoclinic case~{\bf I} will be separate from $1$, then the phase volumes will be unboundedly
contracted or expanded, thus no three-dimensional effective dynamics and Lorenz-like attractors 
will be possible. The same applies to heteroclinic cycles, just here it is reuquired that there exists a  
 point $O_{1}$, where the phase volumes
are contracted, $|J(O_{1})| < 1$ and a point $O_{2}$, near which they are expanded, $|J(O_{2})| > 1$ 
(the so-called, contracting-expanding, 
or mixed case). Controlling the number of iterations the orbit spends near these two points, 
one can avoid unbounded contractions or expansions, i.e. existence of lower-dimensional invariant manifolds.

When both fixed points are saddles (the eigenvalues are real, see Fig.~\ref{saddles}),
one needs to impose some additional bifurcation conditions, 
in order to prevent from existence of lower-dimensional center manifolds. They can be separated into global and local ones.

\subsection{Global degeneracies (non-simple homoclinic and heteroclinic orbits).}
Before giving definitions,
recall some facts from the normal hyperbolicity theory.
Let $O$ be a saddle fixed point of type $(2, 1)$ with eigenvalues $|\lambda_2| < |\lambda_1| < 1 < |\gamma|$ 
and $U_0$ be some small neighbourhood of it.
It is known \cite{GST07, GS90, GS92, book} that
diffeomorphism $\left. f_\mu \right|_{U_0}$ for each small $\mu$ can be represented in some $C^r$-smooth local
coordinates $(x_1, x_2, y)$ as follows (the so-called {\em main normal form}):
\begin{equation}
\begin{array}{l}
\bar x_1 \; = \; \lambda_1(\mu) x_1 +
\tilde H_2(y, \mu)x_2 + O(\|x\|^2|y|)  \\
\bar x_2 \; = \; \lambda_2(\mu) x_2 + \tilde R_2(x, \mu) +
\tilde H_4(y, \mu)x_2 + O(\|x\|^2|y|) \\
\bar y \; = \; \gamma(\mu) y + O(\|x\||y|^2) ,\\
\end{array}
\label{t0norm}
\end{equation}
where $\tilde H_{2,4}(0,\mu) = 0$, $\tilde R_{2}(x, \mu) = O(\|x\|^2)$.
In coordinates (\ref{t0norm}) the invariant manifolds of saddle fixed point $O$ are 
locally straightened: stable $W^s_{loc}(O): \{ y = 0\}$,
unstable $W^u_{loc}(O): \{ x_1 = 0, \; x_2 = 0\}$ and strong stable $W^{ss}_{loc}(O):  \{ x_1 = 0, \; y = 0\}$.

\begin{figure}
\includegraphics[width=8cm]{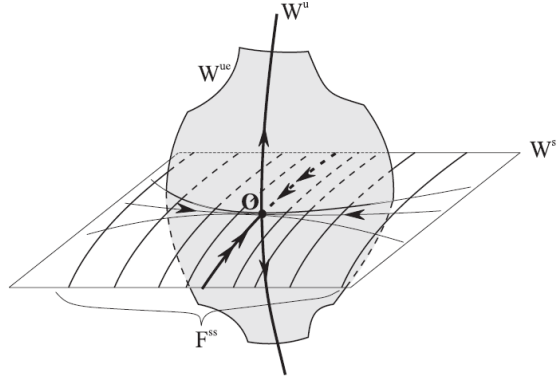}
\caption{\label{fig04} Invariant structures near a saddle fixed point $O$. A part of the strong stable
foliation $F^{ss}$ containing the strong stable manifold $W^{ss}$ and a piece of one of the extended unstable manifolds $W^{ue}$ containing $W^u$ and being transversal to $W^{ss}$ at $O$.} 
\end{figure}

According to Refs.~\onlinecite{book, HPS}, an important role in dynamics is played by an 
{\em extended unstable manifold}  $W^{ue}(O)$,
see Fig.~\ref{fig04}. By definition, it is a two-dimensional invariant manifold, 
that is tangent to the leading stable direction (corresponding to $\lambda_1$) at the saddle point and 
contains unstable manifold $W^u(O)$. Unlike the previous
ones, the extended unstable manifold is not uniquely defined and its smoothness is, generally speaking, 
only $C^{1 + \varepsilon}$. Locally, $W^{ue}_{loc}(O) = W^{ue}(O) \cap U_0$, and the
equation  of $W^{ue}_{loc}(O)$ has the form $x_2 = \varphi(x_1, y)$, where $\varphi(0, y) \equiv 0$ and $\varphi'_{x_1}(0, 0) = 0$.
Note that despite the fact that $W^{ue}(O)$ is non-unique, 
all such manifolds have the same tangent plane at each point of $W^u(O)$.

Another essential fact is the existence of the
{\em strong stable invariant foliation}, see Figure~\ref{fig04}.
In $W^s(O)$ there exists a one-dimensional strong stable invariant submanifold $W^{ss}(O)$, which is
$C^r$--smooth and touches at $O$ the eigenvector corresponding to the strong
stable (nonleading) multiplier $\lambda_2$. Stable manifold $W^s(O)$ is foliated near $O$ by
the leaves of invariant foliation $F^{ss}$ which is $C^r$-smooth, unique and contains
$W^{ss}(O)$ as a leaf.

Now consider a pair of saddle fixed points $O_1$ and $O_2$, and orbit $\Gamma_{21}$ in the points of which
manifolds $W^u(O_2)$ and $W^s(O_1)$ have a quadratic tangency (heteroclinic in case~{\bf II}, when 
$O_1 \neq O_2$, and homoclinic if points $O_1$ and $O_2$ coincide, case~{\bf I}).
Each of the points possesses invariant manifolds and foliations described above.
Let $U_{01} \ni O_1$ and $U_{02} \ni O_2$ be 
some small neighbourhoods of the fixed points, $M_1^+ \in W_{loc}^s(O_1) \subset U_{01}$ and 
$M_2^- \in W_{loc}^u(O_2) \subset U_{02}$ be two points of $\Gamma_{21}$
and $\Pi_1^+ \subset U_{01}$ and $\Pi_2^- \subset U_{02}$ their respective neighborhoods. 
Note that there exists some integer $q_N$ such that $M_1^+ = f_0^{q_1}(M_2^-)$.
The global map along $\Gamma_{21}$ is defined as $T_{21}: \Pi_2^- \to \Pi_1^+ = \left. f^{q_1} \right|_{\Pi_2^-}$.

\begin{figure}
\includegraphics[width=15cm]{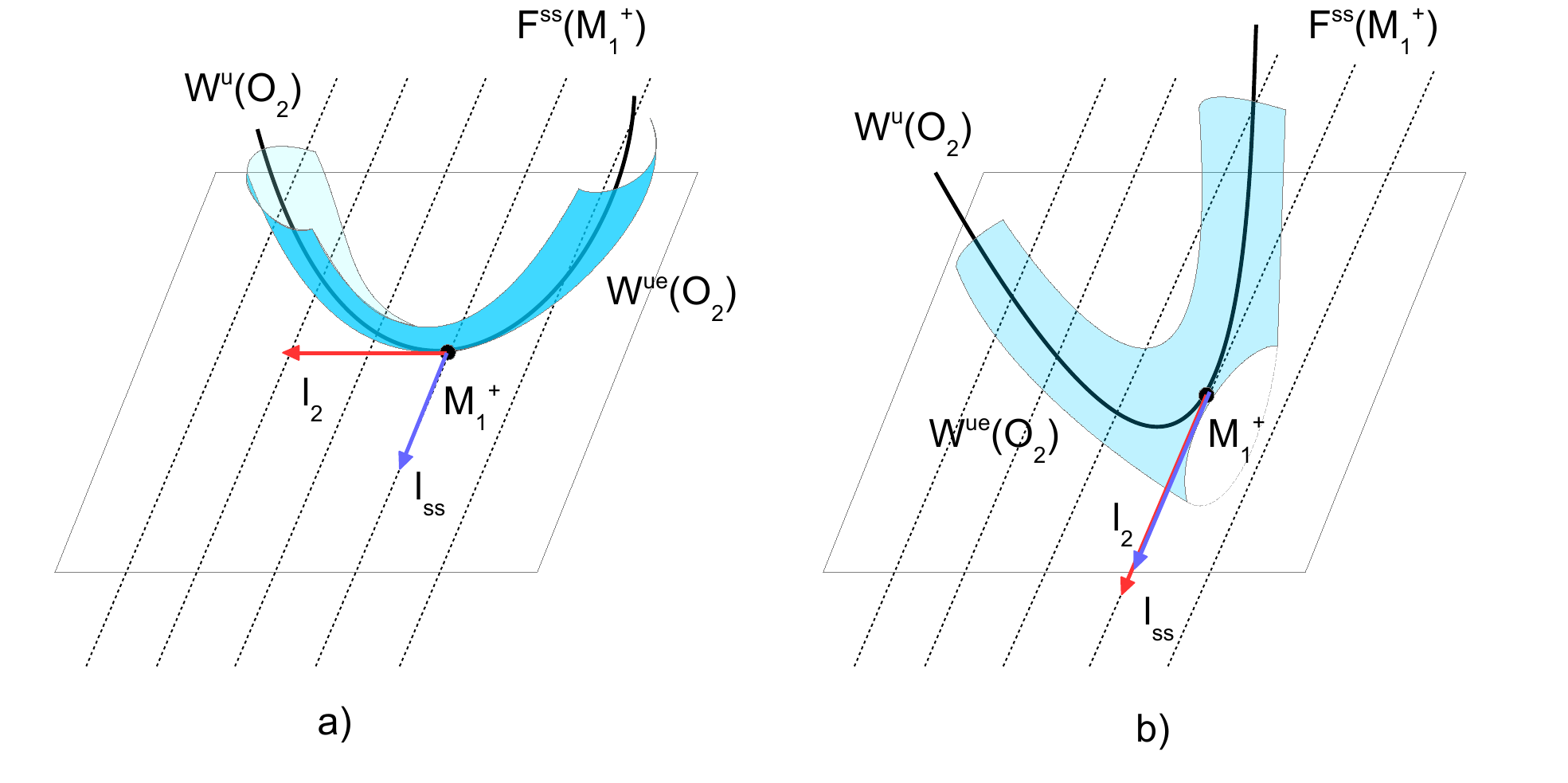}
\caption{\label{fig05} Two types of the non-simple quadratic (homoclinic or heteroclinic) tangency:
(a) Inclination flip: $W^{ue}(O_2)$ is
tangent to $W^s_{loc}(O_1)$ and curves $W^u(O_2)$ and $F^{ss}(M_1^+)$ have a general intersection at
$M_1^+$;
(b) Orbit flip: $W^{ue}(O_2)$ is transversal to $W^s_{loc}(O_1)$ and touches leaf $F^{ss}(M_1^+)$} 
\end{figure}

\begin{df}The homoclinic or heteroclinic tangency of $W^u(O_2)$ and $W^s(O_1)$ is called simple
if image $T_{21}(P^{ue}(M_2^-))$
of tangent plane $P^{ue}(M_2^-)$ to $W^{ue}(O_2)$
intersects transversely the leaf $F_1^{ss}(M_1^+)$ of invariant foliation $F_1^{ss}$, containing point $M_1^+$.
Otherwise, such a quadratic tangency is called non-simple. Following Ref.~\onlinecite{Tat01}, there
may be two generic cases of non-simple homoclinic (heteroclinic) tangencies:

{\rm\bf Orbit flip}. {\it Surface $T_{21}(P^{ue}(M_2^-))$ is
transversal to plane $W^{s}_{loc}(O_1)$  but is tangent to  line
$F_1^{ss}(M_1^+)$ at point $M_1^+$, fig.\ref{fig05}~$($a$)$.}

{\rm\bf Inclination flip}. {\it Surfaces $T_{21}(P^{ue}(M_2^-))$ and
$W^{s}_{loc}(O_1)$ have a quadratic tangency at $M_1^+$ and curves $T_{21}(W^u_{loc}(O_2) \cap \Pi_2^-)$
and $F_1^{ss}(M_1^+)$ have a general intersection, fig.\ref{fig05}~$($b$)$.}  \\
\end{df}

The existence of local extended unstable manifold $W^{ue}_{loc}(O_2)$
implies that near fixed point $O_2$ the dynamics is effectively two-dimensional (the restriction onto
$W^{ue}_{loc}(O_2)$) plus the strong contraction in the transverse direction. If the homoclinic or heteroclinic
tangency is simple, then under forward-time further iterations the image of $T_{21}(P^{ue}(M_2^-))$ will tend
to become tangent to $W^{ue}_{loc}(O_1)$
at point $O_1$. This implies the existence of a global two-dimensional center manifold along $\Gamma_{21}$. When the tangency is non-simple, this global center manifold does not exist.

\begin{remark}
The names for these two degeneracies are taken analogous to the continuous-time case, when the corresponding degeneracies 
also prevent from existence of global two-dimensional center manifolds. The inclination flip corresponds to the case when the extended unstable manifold of $O_2$ has a
quadratic tangency with the stable manifold of $O_1$, similarly to the discrete-time case. 
Orbit flip in flows occurs when the unstable separatrix of $O_2$
comes to $O_1$ along its strong stable direction, that is $W^{ue}(O_2) \supset W^{ss}(O_1)$, and
thus the tangent plane to $W^{ue}(O_2)$ in any point of the homoclinic $($heteroclinic$)$ orbit is also
tangent to $W^{ss}(O_1)$. In discrete-time systems we call it an orbit flip when $W^{ue}(O_2)$ is tangent to any 
leaf of foliation $F^{ss}(O_1)$ (in particular, but not necessary, it can be $W^{ss}(O_1)$ itself). 
Then this tangency will be preserved under the forward iterations of the map, 
therefore images of  $T_{21}(P^{ue}(M_2^-))$ will be always transverse to $W^{ue}_{loc}(O_1)$
and no global center manifold will exist. 
\end{remark}

In the similar way, one defines orbit flip for a transversal heteroclinic orbit $\Gamma_{12}$
(the non-simple heteroclinic intersection).
Consider two points $M_1^- \in U_{01}$ and $M_2^+ \in U_{02}$, 
of  $\Gamma_{12}$ and their small respective neighbourhoods
$\Pi_1^- \subset U_{01}$ and $\Pi_2^+ \subset U_{02}$.
Again, there exists some integer $q_2$ such that $M_2^+ = f_0^{q_2}(M_1^-)$ so that
the global map from $U_{01}$ to $U_2$ is defined as 
$T_{12}: \Pi_1^- \to \Pi_2^+ = \left. f_\mu^{q_2} \right|_{\Pi_1^-}$.
Let $P^{ue}(M_1^-)$ be the tangent plane to $W^{ue}(O_1)$ at $M_1^-$ and $F_2^{ss}(M_{2}^+)$ 
be the leaf of invariant foliation $F_{2}^{ss}$ on $W^s(O_{2})$ passing through $M_{2}^+$. 

\begin{figure}
\includegraphics[width=12cm]{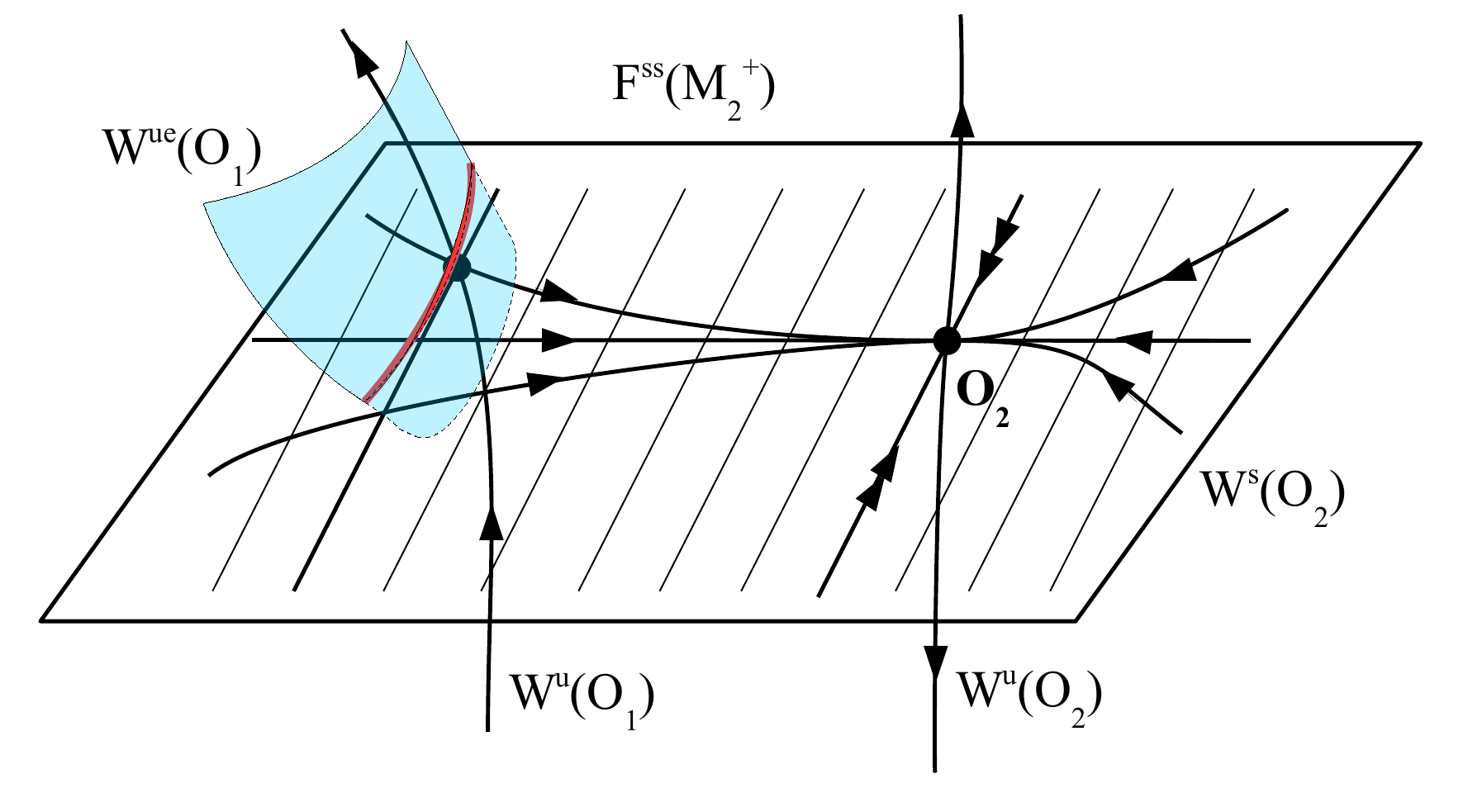}
\caption{\label{Case3} A non-simple heteroclinic intersection (orbit flip) of $W^u(O_1)$ and 
$W^s(O_{2})$.} 
\end{figure}

\begin{df}The heteroclinic intersection of
$W^u(O_1)$ and $W^s(O_{2})$ is called  simple if image $T_{12}(P^{ue}(M_1^-))$ 
and leaf $F_{2}^{ss}(M_{2}^+)$ 
intersect transversely. If this condition is not fulfilled the heteroclinic intersection is non-simple, and undergoes an orbit flip, see fig.~\ref{Case3}.
\end{df} 
%

\subsection{ Local degeneracies.} 
Now consider the case when all fixed points in the heteroclinic cycle are saddles and all connections are simple. 
The existence of extended unstable manifold $W^{ue}$ and strong stable manifold $W^{ss}$ is a robust property --
they persists under small parametric perturbation, and with an absence of non-simple global orbits this immediately implies the existence of a global lower-dimensional center manifold along the homoclinic (heteroclinic) cycle. To prevent from this, local bifurcations can occur at one of the saddle points, such that 
$W^{ue}$ and $W^{ss}$ near this point do not exist at all, cases {\bf I--II.3}. The first
bifurcation is the resonance condition when the stable eigenvalues have the same absolute value, but different signs, 
$\lambda_1 = -\lambda_2 = \lambda$, cases {\bf I--II.3.a}. Under small perturbations here the
strong stable manifold appears, but in alternating directions: when $|\lambda_1| < |\lambda_2|$,
$W^{ss}$ is tangent to the eigendirection corresponding to $\lambda_1$, and to the eigendirection
corresponding to $\lambda_2$ otherwise.

\begin{figure}
\includegraphics[width=10cm]{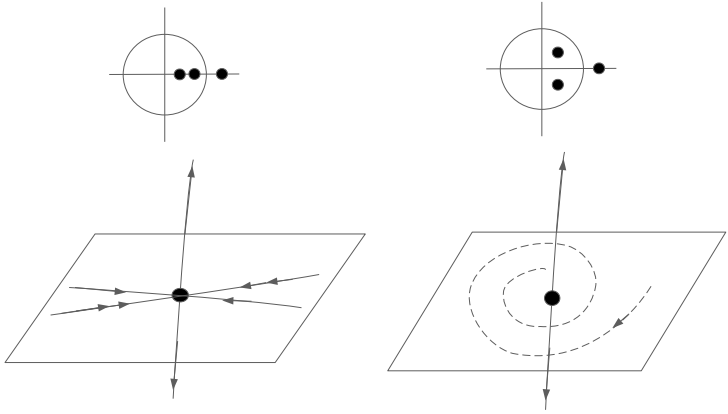}
\caption{\label{fig:transition1} A Belyakov-like transition. When $\mu_2 > 0$ (left), the stable eigenvalues are real, and when $\mu_2 < 0$ (right) they form a complex-conjugate pair.} \end{figure}

The second bifurcation is analogous to the Belyakov resonance for the continuous-time case \cite{Bel80}, that is the boundary between saddle and saddle-focus, cases {\bf I-II.3.b}. 
At the bifurcation moment
the stable multiplier has multiplicity two: $\lambda_1 = \lambda_2 = \lambda$, and under small perturbations
such a degenerate saddle becomes a saddle with real eigenvalues or a saddle-focus, see fig.~\ref{fig:transition1}.

\begin{remark}
Another type of saddle to saddle-focus transition is possible in dimensions higher than three, when a complex-conjugate
pair of stable eigenvalues $\lambda_2 e^{\pm i \varphi}$ coincides in absolute value with real
stable eigenvalue $\lambda_1$, i.e. $|\lambda_1| = \lambda_2 = \lambda$.
In continuous time such a homoclinic bifurcation was studied recently in Ref.~\onlinecite{KMK19}, where it
was called the 3DL-bifurcation. 
Under small perturbations, manifold $W^{ss}$ appears with alternating direction and dimension, 
namely in the case when
$|\lambda_1| < \lambda_2$ it is one-dimensional and tangent to the $\lambda_1$ eigendirection,
and when $|\lambda_1| > \lambda_2$, strong stable manifold $W^{ss}$ is two-dimensional and tangent to 
the eigendirections corresponding to complex eigenvalues $\lambda_2 e^{\pm i \varphi}$. This case is out of scope of the current paper and will be studied separately.
\end{remark}

\subsection{Main theorems}

The main result of the paper is given by the following 
\begin{theorem}  \label{thmmain} 
Let $f_0$ satisfy conditions {\bf A}--{\bf D} and
$f_\mu$ be the three-parametric unfolding  of $f_0$ defined above. Then,
in any neighbourhood of the origin  in the parameter space there exist infinitely many accumulating to $\mu=0$ domains
$\delta_{k}, \; k = \{k_1, k_2\}$, such that diffeomorphism $f_\mu$ has a discrete Lorenz attractor for $\mu \in \delta_{k}$. \\
\end{theorem}

The proof of Theorem~\ref{thmmain} is based on the fact that for sufficiently large $k$ the first return map can be transformed to the form of three-dimensional Henon map~(\ref{H3D}) plus asymptotically small terms, see Lemma~\ref{resc_lemma} for more details. According to Ref.~\onlinecite{GOST05}, map~(\ref{H3D}) possesses the discrete Lorenz attractor is some open domain $V$ in the parameter space. Varying indices $ k_1, k_2$ unboundedly, one gets that in the space of original parameters this domains correspond to a sequence of domains $V_k$  that accumulate to point $\mu = 0$.

\begin{theorem}\label{conj_inv}
Let $f_0$ satisfy conditions {\bf A}--{\bf D} and
$f^{-1}_\mu$ be the three-parametric unfolding  of its inverse map $f^{-1}_0$. Then,
in any neighbourhood of the origin  in the parameter space there exist infinitely many accumulating to $\mu=0$ domains
$\delta_{k}, \; k = \{k_1, k_2\}$, such that diffeomorphism $f^{-1}_\mu$ has a discrete Lorenz attractor for $\mu \in \delta_{k}$. \\
\end{theorem}



\section{Three-dimensional Henon maps and discrete Lorenz attractors}\label{sec:Henon}

In this section the main attention is paid to the extension of results obtained in Refs.~\onlinecite{GOST05} and~\onlinecite{GGOT13} on the birth of discrete Lorenz attractors in local bifurcations of three-dimensional maps. In these papers an analysis of a certain local codimension-three bifurcation was performed, namely when a fixed point possesses multipliers $(-1, -1, +1)$. The reason why 
this bifurcation is relevant to Lorenz attractors
is that the pair of $-1$ multipliers in a single Jordan block provides the same symmetry $(x, y) \to (-x, -y)$ as in the Lorenz system. Also, bifurcations of triply degenerate equilibria with the same symmetry lead to birth of Lorenz attractors in continuous-time systems~\cite{SST93}.

In 3D Henon map (\ref{H3D}), a fixed point with eigenvalues $(-1, -1, +1)$ exists when $(M_1, M_2, B) = (-1/4, 1, 1)$. In paper~\cite{GOST05} the normal form of this bifurcation was approximated by an ODE system, such that for some $T$ a time-$T$ shift of the flow of that system coincided with map (\ref{H3D}) up to arbitrary small terms. It was also shown that this system coincides with the Shimizu-Morioka system after rescaling of coordinates, parameters and time:
\begin{equation}  \label{sm}
      \dot X=Y, \;\; \dot Y=X(1-Z)-\lambda Y, \;\;\dot Z=-\alpha Z+X^2.
\end{equation}
The latter has the Lorenz attractor in some open region of parameters, as proved with assistance of computer, using interval arythmetics in Ref.~\onlinecite{CTZ18}. Thus, the original Henon map can be regarded as the time-$T$ shift of a periodically perturbed ODE system with a Lorenz attractor. From Refs.~\onlinecite{GOST05, GGOT13} it follows that map (\ref{H3D}) possesses a discrete-time analogue of the classical Lorenz attractor.

The results above are applicable to orientation preserving maps, as in the proof the Jacobian $B$ varies near $+1$. When the original map $f_0$ is non-orientable, it is possible that the first return map will also reverse orientation, and parameter $B$ will take only negative values. Then the results above will not apply to such cases, and in order to find Lorenz-like attractors, one needs to look for another parameter domain. Recently, this case was studied in Ref.~\onlinecite{GKKST}. The authors considered a fixed point with eigenvalues $(i, -i, -1)$, which exists in map (\ref{H3D}) for $(M_1, M_2, B) = (7/4, -1, -1)$, and showed that the flow normal form near this bifurcation point possesses a new, 4-winged strange attractor of Lorenz type, which they call the ``Simo's angel''. Numerically such an attractor was observed in Ref.~\onlinecite{GOST05} for $(M_1, M_2, B) = (1.77, -0.925, -0.95)$. Note that the fixed point with eigenvalues $(i, -i, -1)$ has eigenvalues $(-1, -1, +1)$ for the second iterate of the map; however, there is no Jordan block, so the normal form differs significantly from the Shimizu-Morioka system (\ref{sm}). 


In Refs.~\onlinecite{GGOT13, OT17} the results of Ref.~\onlinecite{GOST05} were extended to a wider class of maps and a more generic criterion of existence. It was shown that near a fixed point with multipliers $(-1, -1, +1)$ the map can be represented as the following normal form:
\begin{equation}\label{Lor_nf}
\begin{array}{rcl}
\bar u_1 & = & -u_1 - u_2 \\
\bar u_2 & = & - u_2 + a u_1 u_3 + a_1 u_2 u_3 + O(\|u\|^3) \\
\bar u_3 & = & u_3 + b u_1^2 + b_1 u_2^2 + b_2 u_1 u_2 + b_3 u_3^2 + O(\|u\|^3). \\
\end{array}
\end{equation}
By Lemma 3.1 from Ref.~\onlinecite{GGOT13}, if $ab > 0$, the discrete Lorenz attractors is born in generic perturbations of system (\ref{Lor_nf}), as in this case the flow approximation of it is the Shimizu-Morioka system. From the proof of the Lemma it can be easily seen that when $ab < 0$, the flow normal form can be also transformed to the  Shimizu-Morioka system, but with negative scaling of time. This means that in the normal form a discrete Lorenz repeller is born in generic perturbations. Also it is wotrh to mention, that Lemma 3.1 provides a simple criterion of existence of Lorenz attractors or repellers -- this fact follows immediately from the signs of two coefficients $a$ and $b$ of normal form (\ref{Lor_nf}). This will be used below for the study towards the proof of Conjecture~\ref{conj_inv}.

Now consider the inverse map $f_0^{-1}$, i.e. the diffeomorphism having a homoclinic or heteroclinic cycle composed of $(1, 2)$ saddle or saddle-focus fixed points and one quadratic tangency of invariant manifolds. The first return map $X \to F(X)$ along the cycle after an appropriate rescaling of coordinates and parameters can be brought to the form of an inverse to (\ref{H3D}) map (\ref{Henon2}).
The correspondence between the parameters is 
$\displaystyle \hat B = \frac{1}{B}, \;\;  \hat M_1 = \frac{M_1}{B^2}, \;\; \hat M_2 = -\frac{M_2}{B}$. This map
is also well-known in homoclinic dynamics \cite{GST93c, Tat01, GST96, GST08}. The conservative dynamics of both Henon maps (the case when $B = 1$) was studied in Ref.~\onlinecite{LoMe99}.

As map (\ref{Henon2}) is the inverse to (\ref{H3D}), it automatically follows that it has a discrete Lorenz repeller, which appears under perturbations of a fixed point with multipliers $(-1, -1, +1)$. The ``Simo's angel'' near a fixed point with multipliers $(i, -i, -1)$ is also a repeller here for $\hat B < 0$. It means, that in order to find Lorenz-like attractors in map (\ref{Henon2}), one should look at not only fixed points, but also periodic orbits. That is, to consider $n$-periodic points such that in the $n$-th iterate of map (\ref{Henon2}) it is a fixed point with multipliers $(-1, -1, +1)$ and a Jordan block.

\begin{lm}\label{Lemma_H2}
There exist parameter values, for which map $($\ref{Henon2}$)$ has  periodic points of period $6$ such that its $6$-th iterate $F^6$ has a fixed point with multipliers $(-1, -1, +1)$. The normal form of
this bifurcation is $($\ref{Lor_nf}$)$ with $ab > 0$.
\end{lm}
%
{\it Proof}\\
Theoretical computations show that such periodic points do not exist for periods $2$ and $3$. For periods $4$ and $5$, numerical computations show that there exist parameter values such that $4$th and $5$th iterates of map (\ref{Henon2}) have fixed points with eigenvalues $(-1, -1, +1)$, but in all of them, the coefficients of normal form (\ref{Lor_nf}) give $ab < 0$, this means that discrete Lorenz repellers appear near these points. Now consider period-$6$ orbits. They it will be orbits consisting of
$6$ points $Z_1$--$Z_6$ with coordinates
$$
(z_1, z_2, z_3) \to (z_2, z_3, z_4) \to (z_3, z_4, z_5) \to (z_4, z_5, z_6) \to (z_5, z_6, z_1) \to (z_6, z_1, z_2) \to (z_1, z_2, z_3),
$$ 
which satisfy the following system of equations:
\begin{equation}\label{period-6}
\begin{array}{l}
z_1 = \hat M_1 + \hat M_2 z_6 + \hat B z_4 - z_5^2,  \qquad\qquad
z_2 = \hat M_1 + \hat M_2 z_1 + \hat B z_5 - z_6^2 \\
z_3 = \hat M_1 + \hat M_2 z_2 + \hat B z_6 - z_1^2, \qquad\qquad
z_4 = \hat M_1 + \hat M_2 z_3 + \hat B z_1 - z_2^2 \\
z_5 = \hat M_1 + \hat M_2 z_4 + \hat B z_2 - z_3^2, \qquad\qquad
z_6 = \hat M_1 + \hat M_2 z_5 + \hat B z_3 - z_4^2. \\
\end{array}
\end{equation}

Each of $6$-periodic points $Z_1$--$Z_6$ are fixed points of the sixth iteration of map $F$, i.e. $x \to F^6(x)$. 
Next we will determine the equations that guarantee that matrix $D(F^6)$ possesses eigenvalues $(-1, -1, +1)$ at these points. Namely, condition
\begin{equation}\label{Jacobian}
\det D(F^6)(Z_1) \equiv \hat B^6 = 1
\end{equation}
ensures that the product of the eigenvalues is equal to $1$, condition
\begin{equation}\label{trace}
\tr D(F^6)(Z_1) \equiv \tr \left(DF(Z_6) \circ DF(Z_5) \circ DF(Z_4) \circ DF(Z_3) \circ DF(Z_2) \circ DF(Z_1) \right) = -1 
\end{equation}
makes the sum of the eigenvalues to be equal to $-1$, and the third one
\begin{equation}\label{plus_one}
\det(D(F^6)(Z_1) - \id) = 0
\end{equation}
means that $D(F^6)(Z_1)$ has an eigenvalue $+1$. Then equations (\ref{Jacobian}) and (\ref{trace}) imply that the product of the rest two multiplies is $1$, and their sum is $-2$, so they both are equal to $-1$.

Formulas (\ref{period-6})--(\ref{plus_one}) define $9$ equations for $9$ unknowns $z_1, \ldots, z_6, \hat M_1, \hat M_2, \hat B$, so the system of equations is well-posed. One of numerical solutions of this system is
\begin{equation}\label{p6point-}
\begin{array}{cc}
z_1 = 1.1109087187819051, & z_2 = 0.5430803496704105, \\ 
z_3 = -0.018564282101437988, & z_4 = -1.0126053862814206, \\ 
z_5 = -0.3759675295870319, & z_6 = -0.6947447970072144, \\ 
\hat M_1 = 0.3974562084897318, & \hat M_2 = 0.2271356235631268, \\ \hat B = -1.
\end{array}
\end{equation}

For the parameter values given by (\ref{p6point-}), the $6$-th iterate of the map near point $(z_1, z_2, z_3)$ can be written as normal form (\ref{Lor_nf})
with $a = -0.0555732$ and $b = -1.6955$. According to Ref.~\onlinecite{GGOT13}, Lemma 3.1, a discrete Lorenz attractor (of period $6$) is born in system (\ref{Henon2}) near this bifurcation point for the orientation reversing case, as $\hat B = -1$. 

For the orientation preserving map (\ref{Henon2}), i.e. when $\hat B > 0$, another numerical solution of equations (\ref{period-6})--(\ref{plus_one}) was found:
\begin{equation}\label{p6point+}
\begin{array}{cc}
z_1 =0.913442745966901, & z_2 =1.220643948207064, \\ 
z_3 = 1.3256709760748737, & z_4 = 1.1287783775951246, \\ 
z_5 = 0.7765991221464961, & z_6 = 0.6638157026635255, \\ 
\hat M_1 = -0.9336687216264129, & \hat M_2 =1.99067193080051, \\ \hat B = 1.
\end{array}
\end{equation}
Normal form (\ref{Lor_nf}) has in this case coefficients $a = -0.107789$ and $b = -0.769823$, thus $ab > 0$, and near this bifurcation point also a period-$6$ discrete Lorenz attractor is born.
\mbox



\section{The first return map and the rescaling lemma}\label{sec:fret}

Consider $U$ --- a sufficiently small fixed neighborhood of the homoclinic 
or heteroclinic cycle under consideration. It is a union of small neighbourhoods 
${\textbf U_0} = U_{01}\cup U_{02}$
of the fixed points and small neighbourhoods $U_m$ of all points
of homoclinic or heteroclinic orbits $\Gamma_{12}$ and $\Gamma_{21}$, 
that do not belong 
to ${\textbf U_0}$. Note that there exists only a finite number of such points and neighbourhoods $U_m$. 
Each single-round periodic orbit that lies entirely in
$U$, has exactly one intersection point with each
of $U_m$ and all remaining its points lie in ${\textbf U_0}$.

For each saddle $O_j$, $j = 1,2$, select two points: $M_j^+ \in W^s_{loc}(O_i)$ and $M_j^- \in W^u_{loc}(O_i)$
and their respective neighbourhoods $\Pi_j^+, \Pi_j^- \subset U_{0j}$. The restriction of diffeomorphism
$f_\mu$ onto neighbourhoods $U_{0j}$ are called local maps $T_{0j}$.
Begin iterating $\Pi_j^+$ under the action of $T_{0j}$. Starting from some number $\bar k_j$ images
$T_{0j}^{k} \Pi_j^+$, $k > \bar k_j$, will have nonempty intersections with $\Pi_j^-$. 
As discussed in section~\ref{sec:def}, 
there exist  numbers $q_{1,2}$ such that $M_2^+ = f_0^{q_2}(M_1^-)$, $M_1^+ = f_0^{q_1}(M_2^-)$. For all small $\mu$
global maps are defined as $T_{12}: \Pi_1^- \to \Pi_2^+ = \left. f_\mu^{q_2} \right|_{\Pi_2^-}$, $T_{21}: \Pi_2^- \to \Pi_1^+ = \left. f_\mu^{q_1} \right|_{\Pi_1^-}$.
Now for every $k = (k_1, k_2)$, where $k_j > \bar k_j$, $j = 1,2$,
the first return maps  $T_k: V_k \to \Pi_1^+$ are defined
as $T_k = T_{21} \circ T_{02}^{k_2} \circ T_{12} \circ T_{01}^{k_1}$,
where $V_k \subset \Pi_1^+$ is a subdomain such that 
$T_{01}^{k_1} (V_k) \subset \Pi_1^-$, $T_{12} \circ T_{01}^{k_1} (V_k) \subset \Pi_2^+$, 
$T_{02}^{k_2} \circ T_{12} \circ T_{01}^{k_1} (V_k) \subset \Pi_2^-$, and  $T_k(V_k) \subset \Pi_1^+$.

In order to write the first return map in coordinates, the local and global maps should be represented in the most suitable form. 

\subsection {Local maps}
\begin{table}
\caption{\label{table_local} Local maps near fixed points of different types}
\begin{ruledtabular}
\begin{tabular}{{| c | l | l | c |}}
{\bf NN} & {\bf Fixed point} & {\bf Case \#} & {\bf The local map}  \\
\hline
1. & Saddle & {\bf I.2, II} & (\ref{eq:T0k_s}) \\
\hline
2. & Saddle-focus  & {\bf I--II.1} & (\ref{eq:T0k_sf})\\
\hline
3. & Resonant alternating saddle & {\bf I--II.3.a} & (\ref{eq:T0k_sr})\\
\hline
4. & Resonant Belyakov saddle & {\bf I--II.3.b}& (\ref{eq:T0k_bel})\\
\end{tabular}
\end{ruledtabular}
\end{table}

In this subsection formulas for multiple iterations of local maps, $T_{0j}^{k_j}$ are derived for different types of fixed points:  saddle, saddle-focus, saddle with the alternating resonance, saddle with the Belyakov resonance. The summary with references to formulas is given in Table~\ref{table_local}.

For a saddle fixed point $O_j$ with eigenvalues 
$\lambda_{1j}$, $\lambda_{2j}$, $\gamma_{j}$, $|\lambda_{2j}| < |\lambda_{1j}|$, the local map can be brought to
the main normal form (\ref{t0norm}). This gives us the following formula for its $k$-th iteration
(see Refs.~\onlinecite{book, GS92} for details):
\begin{equation}\label{eq:T0k_s}
\begin{array}{l}
  x_{1k} = \lambda_{(j)1}^k x_{10} +
            \hat\lambda_j^{k} \xi^j_{1k}(x_0, y_k, \mu),\\
   x_{2k} = \hat\lambda_j^{k} \xi^j_{2k}(x_0, y_k, \mu),\\
  y_0  = \gamma_{(j)}^{-k} y_k + \hat \gamma_j^{-k} \xi^j_{3k}(x_0, y_k, \mu).
\end{array}
\end{equation}
Here $0 < |\lambda_{(j)2}| \le \hat \lambda_j < |\lambda_{(j)1}|$, $\hat\gamma_j > |\gamma_{(j)}|$, 
functions $\xi^j_{mk}$ and their derivatives up to the order $(r - 2)$
are uniformly bounded, and their higher order derivatives tend to zero.

Case {\bf I--II.1.} When $O_j$ is a saddle-focus with eigenvalues $\lambda_{(j)} e^{\pm i \varphi_j}$, $\gamma_{(j)}$, 
where $i^2 = -1$, the $k$-th iteration of the local map has the form 
\begin{equation}\label{eq:T0k_sf}
\begin{array}{l}
  (x_{1k}, x_{2k})^\top = \lambda_{(j)}^k R_{k \varphi_j}(x_{10}, x_{20})^\top +
            \hat\lambda_{j}^{k} \xi^j_{1k}(x_0, y_k, \mu),\\
  y_0  = \gamma_{(j)}^{-k} y_k + \hat \gamma_j^{-k} \xi^j_{2k}(x_0, y_k, \mu) ,
\end{array}
\end{equation}
where $R_\psi$ is the rotation matrix of angle $\psi$.

{\bf I--II.3.a.} For the case of a resonant saddle with eigenvalues 
$\lambda_{(j)1}(0) = -\lambda_{(j)2}(0) = \lambda_{(j)}$ and $\gamma_{(j)}$,
the $k$-th iteration can be written as
\begin{equation}\label{eq:T0k_sr}
\begin{array}{l}
  x_{1k} = \lambda_{(j)1}^k x_{10} + \hat\lambda_j^{k} \xi^j_{1k}(x_0, y_k, \mu)\\
   x_{2k} = \lambda_{(j)2}^k x_{20} +\hat\lambda_j^{k} \xi^j_{2k}(x_0, y_k, \mu)\\
  y_0  = \gamma_{(j)}^{-k} y_k + \hat \gamma_j^{-k} \xi^j_{3k}(x_0, y_k, \mu),
\end{array}
\end{equation}
where  $0 < \hat \lambda_j < |\lambda_{(j)}|$. Parameter $\mu_2$ 
unfolds the resonance condition:
\begin{equation}\label{mu2_resonant}
\frac{\lambda_{(j)1}}{\lambda_{(j)2}} = -1 + \mu_2.
\end{equation}

{\bf I--II.3.b.} In the case of the Belyakov-type bifurcation $\lambda_{(j)1}(0) = \lambda_{(j)2}(0) = \lambda_{(j)}$, in order to construct smooth parmetric families,
it is not possible to use canonical Jordan forms for saddle
and saddle-focus~\cite{Arnold_book}, as  
these two normal forms can not be smoothly conjugated at the bifurcation moment. 
One of the possible smooth conjugating parametric families is given by the following formula: 

%
\begin{equation}\label{lin_bel}
Df_\mu(O_j) = \left(\begin{array}{cc}
A_s & 0 \\ 0 & \gamma_{(j)}(\mu)
\end{array}
\right), \;\; 
A_s = \left(\begin{array}{cc}
\lambda_{(j)}(\mu) & 1 \\ \mu_2 & \lambda_{(j)}(\mu) 
\end{array}
\right).
\end{equation}
When $\mu_2 > 0$, the linearization matrix has real stable eigenvalues
$\lambda_{(j)}(\mu) \pm \sqrt{\mu_2}$, and when $\mu_2 < 0$, they form a complex-conjugate pair
$\lambda_{(j)}(\mu) \pm i\sqrt{-\mu_2}$. 
The $k$-th power of matrix $A_s$ has the following form:
\begin{equation}\label{A_k}
A_s^k = \lambda_{(j)}(\mu)^k \left(1 - \frac{\mu_2}{\lambda_{(j)}^2(\mu)}\right)^{k/2} \left(\begin{array}{cc}
C_k(\mu) & S_k(\mu) \\ \mu_2 S_k(\mu) & C_k(\mu)
\end{array}
\right),
\end{equation}
where
\begin{equation}\label{C_k}
C_k = \left\{ \begin{array}{l}
\cosh{k \varphi_j}, \; \mu_2 \ge 0 \\ \cos{k \varphi_j}, \; \mu_2 < 0
\end{array}
\right., \;\;
S_k = \left\{ \begin{array}{l}
\frac{\sinh{k \varphi_j}}{\sqrt{\mu_2}}, \; \mu_2 > 0 \\ 0, \; \mu_2 = 0 \\ 
-\frac{\sin{k \varphi_j}}{\sqrt{-\mu_2}}, \; \mu_2 < 0,
\end{array}
\right. 
\end{equation}
with
\begin{equation}\label{phi_j}
\varphi_j = \left\{ \begin{array}{l}
\arctanh{\frac{\sqrt{\mu_2}}{\lambda_j}}, \; \mu_2 \ge 0 \\ -\arctan{\frac{\sqrt{-\mu_2}}{\lambda_j}}, \; \mu_2 < 0
\end{array}
\right.,
\end{equation}
and the $k$-th iteration of the local map is written as:
\begin{equation}\label{eq:T0k_bel}
\begin{array}{l}
  (x_{1k}, x_{2k})^\top = A_s^k(x_{10}, x_{20})^\top +
            \hat\lambda_j^{k} \xi^j_{1k}(x_0, y_k, \mu),\\
  y_0  = \gamma_{(j)}^{-k} y_k + \hat \gamma_j^{-k} \xi^j_{2k}(x_0, y_k, \mu),
\end{array}
\end{equation}
where $\hat\lambda_j < \lambda_{(j)}$.

All possible types of fixed points and the references to the corresponding formulas for the local maps are given in Table~\ref{table_local}.

\subsection{Global maps}

Recall that the global map in the homoclinic cases maps neighbourhood $\Pi^-$ to $\Pi^+$, and 
in the heteroclinic cases they map $\Pi_1^-$ to $\Pi_2^+$ and $\Pi_2^-$ to $\Pi_{1}^+$.
Assume that the homoclinic or heteroclinic points at $\mu = 0$ have the following coordinates: 
$M_j^-(0, 0, y_{(j)}^-)$ and $M_{j}^+(x_{(j)1}^+, x_{(j)2}^+, 0)$, where $x_{(j)1}^+$, $x_{(j)2}^+$ and $y_{(j)}^-$
depend on parameters, and $(x_{(j)1}^+)^2 + (x_{(j)2}^+)^2 \neq 0$, $y_{(j)}^- \neq 0$. 
The global maps are written as Taylor expansions near points $M_j^+$.

\subsubsection{Transversal heteroclinic intersections}

\begin{table}
\caption{ \label{table_global_transverse} Global maps for transversal intersections}
\begin{ruledtabular}
\begin{tabular}{| c | p{17em} | m{3.8em} | c |}
{\bf NN} &{\bf Connection} & {\bf Case \#} 
& {\bf Genericity conditions} 
\\
\hline
1. & Saddle $\to$ saddle, no orbit flip & {\bf II} &
(\ref{no_orbit_flip_trans})
\\
\hline
2. & Saddle $\to$ saddle, orbit flip & {\bf II.2.c} &
(\ref{mu2_orbit_flip_trans}) 
\\
\hline
3. & Saddle $\to$ saddle-focus \newline 
Saddle-focus $\to$ saddle \newline
Saddle-focus $\to$ saddle-focus & {\bf II.1} & None 
\\
\hline
4. & Resonant alternating saddle $\to$ saddle & {\bf II.3.a} &
(\ref{no_orbit_flip_trans_res_a}) 
 \\
\hline
5. & Saddle $\to$ resonant alternating saddle & {\bf II.3.a} &
(\ref{no_orbit_flip_trans_res_b}) 
\\
\hline

6. & Saddle $\to$ resonant Belyakov saddle & {\bf II.3.b} &
(\ref{no_orbit_flip_trans_resb}) 
\\
\hline
7. & Resonant Belyakov saddle $\to$ saddle & {\bf II.3.b} &
(\ref{no_orbit_flip_trans}) 
 \\
\end{tabular}
\end{ruledtabular}
\end{table}

Consider 
first transversal heterolinic intersections along orbit $\Gamma_{12}$ that appear in 
case {\bf II}. Unstable manifold $W^u(O_1)$ in $\Pi_1^-$ has equation
$x_{(1)1} = x_{(1)2} = 0$, and under the action of global map $T_{12}$ it is transformed to a curve that
intersect transversely stable manifold $W^s(O_2)$, which locally in $\Pi_2^+$ has equation
$y_{(2)} = 0$. Thus one can write $T_{12}$ as follows:
\begin{equation}\label{global_transverse}
\begin{array}{rcl}
x_{(2)1} - x_{(2)1}^+ &=& a_{11}^{(1)} x_{(1)1} + a_{12}^{(1)} x_{(1)2} + b_1^{(1)} (y_{(1)} - y_{(1)}^-) + 
O(\|x_{(1)}\|^2 + |y_{(1)} - y_{(1)}^-|^2) \\
x_{(2)2} - x_{(2)2}^+ &=& a_{21}^{(1)} x_{(1)1} + a_{22}^{(1)} x_{(1)2} + b_2^{(1)} (y_{(1)} - y_{(1)}^-) + 
O(\|x_{(1)}\|^2 + |y_{(1)} - y_{(1)}^-|^2) \\
y_{(2)} &=& y^+_{(1)} + c_{1}^{(1)} x_{(1)1} + c_{2}^{(1)} x_{(1)2} + d^{(1)} (y_{(1)} - y_{(1)}^-) + 
O(\|x_{(1)}\|^2 + |y_{(1)} - y_{(1)}^-|^2). \\
\end{array}
\end{equation}
Here all coefficients depend smoothly on parameters, and $y^+_{(1)}(0) = 0$, $d^{(1)}(0) \neq 0$, as the 
intersection is transversal. Map $T_{12}$ is a diffeomorphism, therefore its Jacobian 
$DT_{12}$ at $M^-_1$ is non-degenerate, i.e.
\begin{equation} \label{Jacobian_transverse}
\det DT_{12} = \det  \left(
\begin{array}{ccc}
a_{11}^{(1)} & a_{12}^{(1)} & b_1^{(1)} \\
a_{21}^{(1)} & a_{22}^{(1)} & b_2^{(1)} \\
c_{1}^{(1)} & c_{2}^{(1)} & d^{(1)}
\end{array}
\right) \neq 0.
\end{equation}

In case {\bf II.2.c} the transversal intersection has an additional degeneracy at the bifurcation moment
--- an orbit flip, in coordinates the condition of a simple and non-simple heteroclinic orbit $\Gamma_{12}$ is 
obtained as follows.
The equation of extended unstable manifold $W^{ue}_{loc}(O_1)$
is $x_{(1)2} = 0$, and the leaf $F^{ss}(M^+_{2})$ passing through point $M^+_{2}$
is locally a straight line $\{x_{(2)1} = x^+_{(2)1}$, $y_{(2)} = 0\}$ with direction vector 
$l^{ss} = (0, 1, 0)^\top$.
Tangent plane $P^{ue}(M^-_1)$ has equation $x_{(1)2} = 0$, and
its image under the action of global map $T_{12}$ has at $\mu = 0$ the following parametric equation:
\begin{equation}\label{global_transverse_1}
\begin{array}{rcl}
x_{(2)1} - x_{(2)1}^+ &=& a_{11}^{(1)} x_{(1)1} + b_1^{(1)} (y_{(1)} - y_{(1)}^-) + 
O(\|x_{(1)}\|^2 + |y_{(1)} - y_{(1)}^-|^2) \\
x_{(2)2} - x_{(2)2}^+ &=& a_{21}^{(1)} x_{(1)1} + b_2^{(1)} (y_{(1)} - y_{(1)}^-) + 
O(\|x_{(1)}\|^2 + |y_{(1)} - y_{(1)}^-|^2) \\
y_{(2)} &=& c_{1}^{(1)} x_{(1)1} + d^{(1)} (y_{(1)} - y_{(1)}^-) + O(\|x_{(1)}\|^2 + |y_{(1)} - y_{(1)}^-|^2). \\
\end{array}
\end{equation}
At point $M^+_{2}$ it has two linearly independent tangent vectors $l_1 = (a_{11}^{(1)}, a_{21}^{(1)}, c_1^{(1)})^\top$
and $l_2 = (b_{1}^{(1)}, b_{2}^{(1)}, d^{(1)})^\top$. Curve $F^{ss}(M^+_{2})$ and surface $T_{12}(P^{ue}(M^-_1))$
will be tangent at point $M^+_{2}$ if vectors $l_1$, $l_2$ and $l^{ss}$ are linearly dependent, this
happens when
\begin{equation}\label{transverse_orbit_flip}
\left. A_{11}^{(1)}(\mu)\right|_{\mu = 0} = 
\left. \left(a_{11}^{(1)}(\mu) - \frac{b_1^{(1)}(\mu) c_1^{(1)}(\mu)}{d^{(1)}(\mu)} \right)\right|_{\mu = 0} = 0.
\end{equation}
So in case {\bf II.2.c}, when the heteroclinic orbit connecting saddles $O_1$ and $O_{2}$ is non-simple,
parameter $\mu_2$ is introduced to unfold the orbit flip degeneracy as
\begin{equation}\label{mu2_orbit_flip_trans}
\mu_2 \equiv A_{11}^{(1)}(\mu).
\end{equation}
When transversal heteroclinic orbit $\Gamma_{12}$ is simple,
it should satisfy the non-degeneracy condition
\begin{equation}\label{no_orbit_flip_trans}
A_{11}^{(1)}(0) \neq 0.
\end{equation}
%

If $O_1$ is a saddle with an alternating resonance (case~{\bf II.3.a}), then due to switching of leading and non-leading
directions for small $\mu$,
$\Gamma_{12}$ will be simple if
\begin{equation}\label{no_orbit_flip_trans_res_a}
A_{11}^{(1)}(0) \neq 0, \;\; 
\left. A_{12}^{(1)}(\mu)\right|_{\mu = 0} = 
\left. \left(a_{12}^{(1)}(\mu) - \frac{b_1^{(1)}(\mu) c_2^{(1)}(\mu)}{d^{(1)}(\mu)} \right)\right|_{\mu = 0} \neq 0.
\end{equation}
If $O_2$ is a saddle with an alternating resonance, the genericity conditions are
\begin{equation}\label{no_orbit_flip_trans_res_b}
A_{11}^{(1)}(0) \neq 0, \;\; 
\left. A_{21}^{(1)}(\mu)\right|_{\mu = 0} = 
\left. \left(a_{21}^{(1)}(\mu) - \frac{b_2^{(1)}(\mu) c_1^{(1)}(\mu)}{d^{(1)}(\mu)} \right)\right|_{\mu = 0} \neq 0.
\end{equation}

If $O_1$ is a saddle undergoing the Belyakov transition (case~{\bf II.3.b}), then the leading stable
direction at $O_1$ tends to the $x_{(1)1}$ axis as $\mu_2 \to +0$, so that condition (\ref{no_orbit_flip_trans})
guarantees the absence of orbit flips in small perturbations. 
If $O_{2}$ undergoes the Belyakov transition, then its
non-leading direction tends to the $x_{(2)1}$ axis in the limit $\mu_2 \to +0$.
In this case the heteroclinic orbit will be simple if
\begin{equation}\label{no_orbit_flip_trans_resb}
A_{21}^{(1)}(0) \neq 0.
\end{equation}

All possible cases of transverse intersections together with the references to the corresponding non-degeneracy conditions are summarized in Table~\ref{table_global_transverse}.

\subsubsection{Quadratic homoclinic and heteroclinic tangencies}

\begin{table} 
\caption{\label{table_global_quadratic} Global maps for quadratic tangencies}
\begin{ruledtabular}
\begin{tabular}{| c | p{17em} | m{4em} | c | }
		{\bf NN} &{\bf Connection} & {\bf Case \#} 
		& {\bf Genericity conditions} \\
		\hline
		1. & Saddle $\to$ saddle, simple tangency & {\bf II} &
		(\ref{simple_tan_1})
		\\
		\hline
		2. & Saddle $\to$ saddle, inclination flip & {\bf I--II.2.a} &
		(\ref{incl_flip}), (\ref{mu2_non-simp}) 
		 \\
		\hline
		3. & Saddle $\to$ saddle, orbit flip & {\bf I--II.2.b} &
		(\ref{orbit_flip}), (\ref{mu2_non-simp})  
		 \\
		\hline
		4. & Saddle $\to$ saddle-focus & {\bf II.1} &
		(\ref{simple_tan_1.1}) 
		 \\
		\hline
		5. & Saddle-focus $\to$ saddle & {\bf II.1} & None 
		\\
		\hline
		6. & Saddle-focus $\to$ saddle-focus & {\bf I--II.1} & None 
		\\
		\hline
		7. & Saddle $\to$ resonant alternating saddle \newline 
		Resonant alternating saddle $\to$ saddle \newline 
		Resonant alternating saddle $\to$ itself& {\bf I--II.3.a} &
		(\ref{simple_tan_2}) 
		 \\
		\hline
		8. & Saddle $\to$ resonant Belyakov saddle \newline
		Resonant Belyakov saddle $\to$ itself 
		& {\bf I--II.3.b} &
		(\ref{simple_tan_3}) 
		 \\
		\hline
		9. & Resonant Belyakov saddle $\to$ saddle & {\bf II.3.b} &
		(\ref{simple_tan_1}) 
		 \\
\end{tabular}
\end{ruledtabular}	
\end{table}

The nontransversal heteroclinic (homoclinic) orbit connects 
fixed points $O_2$ and $O_1$. 
When $\mu = 0$, global map $T_{21}$ transforms a piece of unstable manifold 
$W^{u}(O_2) \cap \Pi^-_2$ with equation
$x = 0$ into a curve tangent at point $M^+_1$ to surface $W^s(O_1) \cap \Pi^+_1$ with equation $y = 0$. Then for
all small $\mu$  global map $T_{21}$ is written as 
\begin{equation}\label{global_tangency_3D}
\begin{array}{rcl}
\bar x_{(1)1} - x_{(1)1}^+ &=& a_{11}^2 x_{(2)1} + a_{12}^2 x_{(2)2} + b_1^2 (y_{(2)} - y_{(2)}^-) + 
O(\|x_{(2)}\|^2 + |y_{(2)} - y_{(2)}^-|^2) \\
\bar x_{(1)2} - x_{(1)2}^+ &=& a_{21}^2 x_{(2)1} + a_{22}^2 x_{(2)2} + b_2^2 (y_{(2)} - y_{(2)}^-) + 
O(\|x_{(2)}\|^2 + |y_{(2)} - y_{(2)}^-|^2) \\
\bar y_{(1)1} &=& y^+_{(2)} + c_{1}^2 x_{(2)1} + c_{2}^2 x_{(2)2} + d^2 (y_{(2)} - y_{(2)}^-)^2 + 
O(\|x_{(2)}\|^2 + |y_{(2)} - y_{(2)}^-|^3) \\
\end{array}
\end{equation}

The left hand side variables are denoted here
as $(\bar x_{(1)}, \bar y_{(1)})$ to indicate that the image of $T_{21}$ lies in $\Pi_1^+$ and these coordinates
also represent the iteration of the first return map $T_k$ from $\Pi_1^+$ to itself. All coefficients here depend smoothly on parameters, and when $\mu = 0$ we have $y^+_{(2)}(0) = 0$ and
$d^2(0) \neq 0$, as the tangency is quadratic at the bifurcation moment. The Jacobian of the global map
$DT_{21}$ at $M^-_2$ is non-degenerate, that is
\begin{equation}\label{jacobian_tangency}
\det DT_{21} = \det  \left(
\begin{array}{ccc}
a_{11}^2 & a_{12}^2 & b_1^2 \\
a_{21}^2 & a_{22}^2 & b_2^2 \\
c_{1}^2 & c_{2}^2 & 0
\end{array}
\right) \neq 0
\end{equation}

When $\mu \neq 0$, value $y^+_{(2)}(\mu)$ is the splitting distance of the quadratic tangency up to $o\|\mu\|$
terms, so it is taken as the splitting parameter:
\begin{equation}\label{mu_1}
\mu_1 \equiv y^+_{(2)}(\mu).
\end{equation}

Now write in coordinates the conditions of simple and non-simple quadratic tangencies. Consider saddle fixed points $O_2$ and $O_1$ such that all their eigenvalues are real and do not satisfy resonance conditions from cases {\bf I--II.3}.  
The equation of extended unstable manifold $W^{ue}_{loc}(O_2)$
is $x_{(2)2} = 0$, and the leaf $F^{ss}(M^+_{1})$ passing through point $M^+_{1}$
is locally a straight line $\bar x_{(1)1} = x^+_{(1)1}$, $\bar y_{(1)} = 0$ with direction vector 
$l^{ss} = (0, 1, 0)^\top$.
The image of tangent plane $P^{ue}(M^-_2)$ under the action of global map $T_{21}$ has the following
parametric equation:
\begin{equation}\label{global_tangency_1}
\begin{array}{rcl}
\bar x_{(1)1} - x_{(1)1}^+ &=& a_{11}^2 x_{(2)1} + b_1^2 (y_{(2)} - y_{(2)}^-) + 
O(\|x_{(2)}\|^2 + |y_{(2)} - y_{(2)}^-|^2) \\
\bar x_{(1)2} - x_{(1)2}^+ &=& a_{21}^2 x_{(2)1} + b_2^2 (y_{(2)} - y_{(2)}^-) + 
O(\|x_{(2)}\|^2 + |y_{(2)} - y_{(2)}^-|^2) \\
\bar y_{(1)} &=& c_{1}^2 x_{(2)1} + d^2 (y_{(2)} - y_{(2)}^-)^2 + O(\|x_{(2)}\|^2 + |y_{(2)} - y_{(2)}^-|^2). \\
\end{array}
\end{equation}
At point $M^+_{1}$ it has two linearly independent tangent vectors $l_1 = (a_{11}^2, a_{21}^2, c_1^2)^\top$
and $l_2 = (b_{1}^2, b_{2}^2, 0)^\top$. Curve $F^{ss}(M^+_{1})$ and surface $T_{21}(P^{ue}(M^-_2))$
will be tangent at point $M^+_{1}$ if vectors $l_1$, $l_2$ and $l^{ss}$ are linearly dependent, this
happens when
\begin{equation}\label{tangency_non-simp}
\left. b_1^{2}(\mu) c_1^{2}(\mu)\right|_{\mu = 0} = 0.
\end{equation}
So here naturally two possibilities appear for the quadratic tangency to be non-simple. 
In the inclination flip cases {\bf I--II.2.a}, surfaces $T_{21}(P^{ue}(M^-_2))$ and
$W^s(O_1)$ are tangent to each other (fig.~\ref{fig05}~(a)), therefore vectors $l_1$ and $l_2$ both lie in $W^s(O_1)$, 
thus
\begin{equation}\label{incl_flip}
c_1^2(0) = 0, \; b_1^2(0) \neq 0,
\end{equation}
and in the orbit flip cases {\bf I--II.2.b}, when surface $T_{21}(P^{ue}(M^-_2))$ is transverse
to $W^s(O_1)$, (fig.~\ref{fig05}~(b)), it follows that
\begin{equation}\label{orbit_flip}
b_1^2(0) = 0, \; c_1^2(0) \neq 0,
\end{equation}
which implies that vectors $l_{ss}$ and $l_2$ are parallel.

For these types of degeneracies the second unfolding parameter $\mu_2$ is introduced as
\begin{equation}\label{mu2_non-simp}
\mu_2 = \left\{\begin{array}{rl}
c_1^2(\mu) & {\rm in \; cases \; \textbf{I--II.2.a}} \\
b_1^2(\mu) & {\rm in \; cases \; \textbf{I--II.2.b}} 
\end{array}\right.
\end{equation}

When the quadratic tangency is simple
the condition of absence of non-simple tangencies at the bifurcation moment and in small perturbations should be written. 
If points $O_1$ and $O_2$ are saddles, and they do not satisfy resonance conditions from cases
cases {\bf I--II.3}, then
\begin{equation}\label{simple_tan_1}
b_1^2(0) \neq 0, \; c_1^2(0) \neq 0.
\end{equation}

If point $O_2$ is a saddle and $O_1$ is a saddle-focus, only the inclination flip degeneracy is possible, when manifold $W^{ue}(O_2)$
is tangent to stable manifold $W^s(O_1)$. To avoid this, one needs:
\begin{equation}\label{simple_tan_1.1}
c_1^2(0) \neq 0.
\end{equation}

If one of the points $O_1$ and $O_2$ at the bifurcation moment is a saddle with the alternating resonance,
cases {\bf I--II.3.a}, then either the direction of $W^{ue}(O_2)$, or the direction of $W^{ss}(O_1)$
may alternate when $\mu$ varies, thus the quadratic tangency is simple if
\begin{equation}\label{simple_tan_2}
b_1^2(0) \neq 0, \; b_2^2(0) \neq 0, \; c_1^2(0) \neq 0, \; c_2^2(0) \neq 0.
\end{equation}

If point $O_2$ satisfies the Belyakov condition, case {\bf II.3.b}, then 
inequalities (\ref{simple_tan_1}) should be fulfilled to avoid non-simple tangencies, and if $O_1$ satisfies the Belyakov condition 
(this also inlcludes the homoclinic case {\bf I.3.b}), the quadratic tangency will be simple if
\begin{equation}\label{simple_tan_3}
b_2^2(0) \neq 0, \; c_1^2(0) \neq 0.
\end{equation}

All possible cases of quadratic tangencies together with the references to the corresponding non-degeneracy conditions are summarized in Table~\ref{table_global_quadratic}.

\begin{lm}\label{resc_lemma} {\em (The rescaling lemma)}
Let $f_{\mu_1, \mu_2, \mu_3}$ be the family under consideration. Then, in the
space $(\mu_1, \mu_2, \mu_3)$ there exist infinitely many regions $\Delta_{i}$ in the homoclinic case {\bf I} and $\Delta_{ij}$ in the heteroclinic case {\bf II}
ac\-cu\-mu\-la\-ting to the origin as $i, j \to \infty$, such that the first return map in appropriate rescaled
coordinates and parameters is asymptotically $C^{r - 1}$-close to one of the following limit maps.

{\rm 1)} In the orbit flip cases~{\bf I--II.2.b}:
\begin{equation}
\label{HIhet}
\begin{array}{l}
\bar X_1 \; = \; -B X_2 + M_2 Y,\;\; \bar X_2 \; = \; Y,\;\;
\bar Y = M_1 - X_1 - Y^2,
\end{array}
\end{equation}

{\rm 2)}  In all other cases:
\begin{equation}
\label{HIIhet}
\begin{array}{l}
\bar X_1 \; = \; Y,\;\; \bar X_2 \; = \; X_1,\;\;
\bar Y = M_1 + M_2 X_1 + B X_2- Y^2,
\end{array}
\end{equation}

\end{lm}

Thus, the rescaled first return map in almost all cases is exactly the 3D Henon map (\ref{H3D}). In
cases~{\bf I--II.2.b} in system (\ref{HIhet}) we make an additional change of coordinates $X_{1new} = X_1 - M_2 X_2$ and scale $X_1$ by $(-B)$, bringing it
again to the form (\ref{HIIhet}). 

The relations between old and new parameters are the following.
\begin{equation}\label{mu1_resc}
	M_1 \sim \left\{\begin{array}{rl}
		\mu_1 \gamma^{2i} & {\rm in \; case \; \textbf{I}} \\
		\mu_1 \gamma_{(1)}^{2i} \gamma_{(2)}^{2j} & {\rm in \; case \; \textbf{II}}.
	\end{array}\right.
\end{equation}
When $i, j \to \infty$, with sufficiently small variations of parameter $\mu_1$ one gets arbitrary finite values of parameter $M_1$.

\begin{equation}\label{mu3_resc}
	B \sim \left\{\begin{array}{rl}
		J^i(O) \det DT_1 & {\rm in \; case \; \textbf{I}} \\
		J^i(O_1)J^j(O_2)\det DT_{12} \det DT_{21} & {\rm in \; case \; \textbf{II}}. 
	\end{array}\right.
\end{equation}
Based on formulas (\ref{mu3_hom}) and (\ref{mu3_het}), by small variations of parameter $\mu_3$ parameter $B$ takes arbitrary finite values. If the original diffeomorphism $f_0$ is orientation preserving, $B$ takes only positive values, if $f_0$ is orientation reversing, then $B$ takes either only positive or only negative values, depending on the orientability of the first return map.
\begin{equation}\label{mu2_resc_1}
	M_2 \sim \left\{\begin{array}{rl}
		\lambda_1^i \gamma^i \cos(i \varphi + \theta) & {\rm in \; case \; \textbf{I.1}} \\
		\lambda_{(1)1}^i \gamma_{(1)}^i \lambda_{(2)1}^j \gamma_{(2)}^j \cos(i \varphi_1 + \theta_1) \cos(j \varphi_2 + \theta_2) & {\rm in \; case \; \textbf{II.1}}, 
	\end{array}\right.
\end{equation}
where $\theta$, $\theta_1$ and $\theta_2$ smoothly depend on parameters and $\mu_2$ is varied in the way that the trigonometric function stays close to zero. At the same time, according to formulas (\ref{mu3_hom}) and (\ref{mu3_het}), the coefficients
$$
\lambda_1^i \gamma^i \sim \lambda_2^{-i} \;\;
{\rm and } \;\;
\lambda_{(1)1}^i \gamma_{(1)}^i \lambda_{(2)1}^j \gamma_{(2)}^j \sim \lambda_{(1)2}^{-i} \lambda_{(2)2}^{-j}
$$
are asymptotically large when $i, j \to \infty$. Thus parameter $M_2$ takes arbitrary finite values.
\begin{equation}\label{mu2_resc_2}
	M_2 \sim \left\{\begin{array}{rl}
		\mu_2 \lambda_1^i \gamma^i  & {\rm in \; case \; \textbf{I.2}} \\
		\mu_2 \lambda_{(1)1}^i \gamma_{(1)}^i \lambda_{(2)1}^j \gamma_{(2)}^j  & {\rm in \; case \; \textbf{II.2}}. 
	\end{array}\right.
\end{equation}
Again, for  $i, j \to \infty$ and sufficiently small $\mu_2$ parameter $M_2$ takes arbitrary finite values.
\begin{equation}\label{mu2_resc_3}
	M_2 \sim \left\{\begin{array}{rl}
		\displaystyle
		\lambda_1^i \gamma^i \left(
		(-1 + \mu_2)^i+ A \right) & {\rm in \; case \; \textbf{I.3.a}} \\
		 \displaystyle
		 \lambda_{(1)1}^i \gamma_{(1)}^i \lambda_{(2)1}^j \gamma_{(2)}^j \left(
		 (-1 + \mu_2)^k 
		 + A \right) & {\rm in \; case \; \textbf{II.3.a}}.\\
\end{array}\right.
\end{equation}
Here value $A$  smoothly depends on the parameters, and $A \neq 0$ when $\mu = 0$. The power $k$ denotes $i$ or $j$ depending on which saddle point, $O_1$ or $O_2$, satisfies  the resonance condition. Here, to make $M_2$ finite, parameter $\mu_2$ is varied near such values, where $\left((-1 + \mu_2)^k + A \right)$ becomes zero. To achieve this, the parity of $k$ is taken appropriately, depending on the sign of $A$.
\begin{equation}\label{mu2_resc_4}
	M_2 \sim \left\{\begin{array}{rl}
		\displaystyle
		\lambda_1^i \gamma^i \left(\frac{A}{\sqrt{-\mu_2}}\cos(i \varphi + \theta) \right)  & {\rm in \; case \; \textbf{I.3.b}} \\
		\displaystyle
		\lambda_{(1)1}^i \gamma_{(1)}^i \lambda_{(2)1}^j \gamma_{(2)}^j \left(\frac{A}{\sqrt{-\mu_2}}\cos(k \varphi + \theta) \right) & {\rm in \; case \; \textbf{II.3.b}}. 
	\end{array}\right.
\end{equation}
This formula is valid only when $\mu_2 < 0$, which means that the saddle point having a stable eigenvalue with multiplicity two (the Belyakov resonance), becomes a saddle-focus. Here $A$ and $\theta$ smoothly depend on the parameters, moreover $A \neq 0$ when $\mu = 0$. Exponent $k$ is $i$ or $j$ depending on which saddle point, $O_1$ or $O_2$, satisfies the resonance condition. 
The angle variable $\varphi$ is given by formula (\ref{phi_j}). Varying a small $\mu_2$ near one of the zeros of the trigonometric function, and at the same time, keeping it away from zero, one get parameter $M_2$ taking arbitrary finite values.


\section*{Acknowledgements}
This paper is a contribution to the project M7 (Dynamics of Geophysical Problems in Turbulent Regimes) of the Collaborative Research Centre TRR 181 ``Energy Transfer in Atmosphere and Ocean'' funded by the Deutsche Forschungsgemeinschaft (DFG, German Research Foundation) -- Projektnummer 274762653. The paper is also supported by the
grant of the Russian Science
Foundation 19-11-00280. I especially thank S. Gonchenko and D. Turaev for the idea of considering bifurcations of periodic points in the inverse 3D Henon map (\ref{Henon2}). In addition, I am grateful to Jean-Luc Doumont for the seminars on scientific writing, that helped me to significantly improve the manuscript.

\section*{Data Availability Statement}
The data that support the findings of this study (proofs, pictures) are placed in the body of the
text. If some extra requirements appear, they should be addressed to the corresponding author.

\nocite{*}
\bibliography{Ovsyannikov}

\end{document}